\documentclass[pdflatex,sn-mathphys-num]{sn-jnl}
\usepackage{amsmath,amssymb,amsthm,mathtools}
\usepackage[T1]{fontenc}
\usepackage{booktabs,tikz,placeins}
\usetikzlibrary{arrows.meta}
\usepackage{graphicx}
\hypersetup{colorlinks=true,linkcolor=black,citecolor=black,urlcolor=black}
\makeatletter
\let\arta@headings\ps@headings
\def\ps@headings{\arta@headings
  \def\@oddhead{\hfil{\small\itshape\s@title}\hfil}
  \let\@evenhead\@oddhead}
\def\pacsbullet{\hbox{,\hskip2.5pt}}
\let\arta@titlepage\ps@titlepage
\def\ps@titlepage{\arta@titlepage\let\@oddhead\@empty\let\@evenhead\@empty}
\makeatother
\numberwithin{equation}{section}
\newtheorem{theorem}{Theorem}[section]
\newtheorem{proposition}[theorem]{Proposition}
\newtheorem{lemma}[theorem]{Lemma}
\newtheorem{corollary}[theorem]{Corollary}
\theoremstyle{definition}
\newtheorem{definition}[theorem]{Definition}

\theoremstyle{remark}\newtheorem{remark}[theorem]{Remark}
\newcommand{\mathscrS}{\mathcal{S}_0}
\newcommand{\C}{\mathcal C}
\newcommand{\D}{\mathcal D}
\newcommand{\E}{\mathbb E}
\newcommand{\Hom}{\operatorname{Hom}}
\newcommand{\End}{\operatorname{End}}
\newcommand{\Ext}{\operatorname{Ext}}
\newcommand{\rad}{\operatorname{rad}}
\newcommand{\ind}{\operatorname{ind}}
\newcommand{\add}{\operatorname{add}}
\newcommand{\modu}{\operatorname{mod}}
\newcommand{\Ex}{\operatorname{Ex}}
\newcommand{\AR}{\operatorname{AR}}
\newcommand{\coker}{\operatorname{Coker}}
\newcommand{\eu}{\operatorname{eu}}
\newcommand{\ksp}{K_0^{\mathrm{sp}}}
\newcommand{\dimk}{\dim_k}
\title[Grothendieck groups of extended module categories]{Grothendieck groups and Auslander--Reiten relations in extended module categories}
\newif\ifShowAuthorInformation
\ShowAuthorInformationtrue
\newcommand{\AuthorOneName}{Yajie Xiang}
\newcommand{\AffiliationOne}{School of General Education, Beijing Vocational College of Labour and Social Security, 100029 Beijing, China}
\newif\ifShowSecondAuthor
\ShowSecondAuthortrue
\newcommand{\AuthorTwoName}{Xiao Zhuang}
\newcommand{\AffiliationTwo}{Institute of Fundamental and Interdisciplinary Sciences, Beijing Union University, 100101 Beijing, China}
\newcommand{\CorrespondingAuthorNumber}{2}
\newcommand{\CorrespondingAuthorEmail}{ldtzhuangxiao@buu.edu.cn}
\newcommand{\CorrespondenceDetails}{}
\newcommand{\FundingStatement}{Xiao Zhuang was supported by the NSF of China (Grant No. 12201045), the R\&D Program of Beijing
Municipal Education Commission (Grant No. KM202211417002),the Beijing Union University Science and Technology Innovation
Incubation Program Funding (Grant No. JZK30202506).}
\newcommand{\AIUseStatement}{We acknowledge the use of AI during the preparation of this manuscript. The main proof strategies were developed by the authors. ChatGPT assisted primarily with the counting of indecomposable objects in the extended module categories of radical-square-zero linear Nakayama algebras considered in Section 5.2. The authors have reviewed the manuscript in full and take full responsibility for all mathematical claims, arguments, and conclusions, as well as for any remaining errors.}

\ifShowAuthorInformation
  \ifnum\CorrespondingAuthorNumber=1
    \author*[1]{\AuthorOneName}\email{\CorrespondingAuthorEmail}
  \else
    \author[1]{\AuthorOneName}
  \fi
  \ifShowSecondAuthor
    \ifnum\CorrespondingAuthorNumber=2
      \author*[2]{\AuthorTwoName}\email{\CorrespondingAuthorEmail}
    \else
      \author[2]{\AuthorTwoName}
 \fi
  \affil[1]{\AffiliationOne}
  \ifShowSecondAuthor\affil[2]{\AffiliationTwo}\fi
    
\fi

\newcommand{\PrintDeclaration}[2]{%
  \ifx#2\empty\else\subsection*{#1}#2\par\fi}
\newcommand{\FrontFunding}{%
  \ifx\FundingStatement\empty\else
    \begingroup\small\noindent\textbf{Funding.} \FundingStatement\par\medskip\endgroup
  \fi}
\begin{document}
\raggedbottom
\abstract{
In this paper we consider Grothendieck groups of extended module
categories $\C_m(A)=m\text{-}\modu A$, where $A$ is a finite-dimensional
algebra over a field and $m$ is a positive integer. We show that the
relation subgroup of the split Grothendieck group of $\C_m(A)$ is
generated by Auslander--Reiten extriangles if and only if $\C_m(A)$
has finitely many indecomposable objects up to isomorphism. To this end,
we establish a criterion for extriangulated categories with bounded
negative extension sequences and a finite extension generator. We then
apply this criterion to extended module categories and recover the
two-term silting criterion of Padrol, Palu, Pilaud and Plamondon.
We also use the cohomological Euler characteristic to determine the
rank of the relation subgroup. In particular, when $A$ is hereditary
and representation-finite, we explicitly describe the indecomposable
objects of $\C_m(A)$ and calculate the rank of its relation subgroup.
}
\maketitle
\thispagestyle{titlepage}
\ifShowAuthorInformation
\ifx\CorrespondenceDetails\empty\else
\noindent{\small\textbf{Correspondence.} \CorrespondenceDetails}\par\medskip
\fi
\fi
\FrontFunding

\section{Introduction}
Butler \cite{Butler} and Auslander \cite{Auslander} characterized
representation-finite Artin algebras by the generation of their
Grothendieck relations by almost split sequences. Enomoto
\cite{Enomoto} studied this question for exact categories using effaceable functors
and finiteness of defects. For triangulated categories, Xiao and Zhu
\cite{XZrelations,XZlocal} proved generation results, and Haugland
\cite{Haugland} established converses under
appropriate finiteness hypotheses. Fedele studied Auslander--Reiten $(d+2)$-angles in subcategories
\cite{FedeleAR} and Grothendieck groups of triangulated categories
via cluster tilting subcategories \cite{FedeleK}. Zhou \cite{PZhou}
studied the same question for Hom-finite Krull--Schmidt
$(d+2)$-angulated categories over an algebraically closed field with
Auslander--Reiten $(d+2)$-angles. For odd $d$, he proved that
generation of Grothendieck relations by those angles implies local
finiteness. He, Liu and Zhou
\cite{HLZ} proved the
even-dimensional counterpart under an additional non-cancellation
hypothesis.

Nakaoka and Palu \cite{NP} introduced extriangulated categories, which
provide a common framework for exact and triangulated categories.
Iyama, Nakaoka and Palu \cite{INP} developed Auslander--Reiten theory
in this setting. Zhu and Zhuang \cite{ZZ} proved that,
for a locally finite extriangulated category, the relation subgroup
of its split Grothendieck group is generated by Auslander--Reiten
extriangles (AR extriangles for short). Padrol, Palu, Pilaud and Plamondon
\cite{PPPP}, abbreviated to PPPP,
established a related finiteness criterion and relation-basis theorem
for categories admitting two-term projective presentations.

Our main application concerns the extended module categories studied by
Gupta \cite{Gupta} and Zhou \cite{Zhou}; see also Gupta and
Zhou \cite{GZ}. Let $k$ be an
arbitrary field, let $A$ be a finite-dimensional unital $k$-algebra, and
let $m\geq1$. We use finitely generated right modules and set
\[
 \D=D^b(\modu A),\qquad
 \C_m(A)=\{X\in\D:H^q(X)=0\text{ for }q\notin[1-m,0]\},
\]
with $H^q(X[1])=H^{q+1}(X)$. This extension-closed subcategory inherits
an extriangulated structure. For $m\geq2$ and $A\ne0$, the extriangle
$M\to0\to M[1]$ with $M\ne0$ shows that this structure is not exact;
the category is also not closed under the ambient shift.

Write $\ksp(\C)$ for the split Grothendieck group. An extriangle
$\delta:X\to Y\to Z$ defines $r_\delta=[X]-[Y]+[Z]$; the subgroups
generated by all such relations and by Auslander--Reiten relations are
denoted $\Ex(\C)$ and $\AR(\C)$, respectively. Choose simple modules
$S_1,\ldots,S_n$, put $S=\bigoplus_aS_a$, and choose an injective
embedding $S\hookrightarrow J$ with cokernel $Q$. Direct sums of
shifted short exact sequences and shift triangles give an extriangle
\begin{equation}\label{eq:introdetector}
 \epsilon:F\longrightarrow J[m-1]\longrightarrow T,\qquad
 F=\bigoplus_{i=0}^{m-1}S[i],\quad
 T=\bigoplus_{i=1}^{m-1}S[i]\oplus Q[m-1].
\end{equation}

\begin{theorem}[Theorem~\ref{thm:finite}]\label{intro:finite}
The following conditions are equivalent:
\begin{enumerate}
\item $\C_m(A)$ has finitely many indecomposable objects up to isomorphism;
\item $\C_m(A)$ is locally finite;
\item $\Ex(\C_m(A))=\AR(\C_m(A))$;
\item $r_\epsilon\in\AR(\C_m(A))$.
\end{enumerate}
\end{theorem}

The main step is to associate to every indecomposable object $U$ a
function on the split Grothendieck group satisfying
\[
 \chi^-_U(r_\delta)=\dimk D_\delta(U),\qquad
 \chi^-_U(V)=\sum_{j=0}^{m-1}(-1)^j\dimk\Hom_{\D}(U,V[-j]),
\]
where $D_\delta$ is the right defect of the extriangle $\delta$.
For an AR extriangle ending at $Z$, this value is zero unless
$U\simeq Z$. In contrast, its value on $r_\epsilon$ is positive for
every non-projective indecomposable $U$. If $r_\epsilon$ is a finite
sum of AR relations, all non-projective indecomposables must therefore
occur among the finitely many endpoints in that sum. Since there are
also only finitely many indecomposable projectives, the category is finite.
This use of one relation follows the arguments of Haugland
\cite{Haugland} and Zhu and Zhuang
\cite{ZZ}. PPPP
\cite{PPPP} use the relation $[P]+[P[1]]$
in the two-term setting, while in our general criterion, the required
relation is obtained from an extriangle whose first term is a finite
extension generator and whose middle term is injective.

When the category is finite, the AR relations form an integral basis
of the relation subgroup. The coefficient of the relation ending at
$Z$ is
\[
 \frac{\dimk D_\delta(Z)}{\dimk(\End(Z)/\rad\End(Z))}.
\]
PPPP \cite{PPPP} obtained a basis
expansion using the starting objects of AR extriangles and a quotient
of Hom spaces. The formula above uses their endpoints and right
defects. The bounded negative extension sequence supplies the
dimension identity needed for this formula. We also identify
$K_0(\C_m(A))$ with $K_0(\modu A)$ by the cohomological Euler
characteristic, following the argument of the Stacks Project Authors
\cite[Lemma 13.28.2, Tag 0FCP]{Stacks}, and determine the full relation
subgroup and its rank. After submitting this paper, we learned that Zhou and Zhu\cite{ZhouZhu} had obtained a general characterization of local finiteness for extriangulated categories in terms of Grothendieck groups; while the present paper expresses defect dimensions through negative extensions and also proves that locally finite is equivalent to finite for extended module categories.

The paper is organized as follows. In Section~\ref{sec:preliminaries},
we recall extriangulated categories, their basic exact sequences,
and the notation for extended module categories and Grothendieck
groups. In Section~\ref{sec:general}, we prove the finiteness criterion
for categories with bounded negative extension sequences and derive
the AR relation-basis formula. In Section~\ref{sec:module}, we apply
the criterion to extended module categories, calculate their
Grothendieck groups and relation ranks, and recover the two-term
silting criterion. In Section~\ref{sec:examples}, we compute a
non-hereditary Nakayama example, count indecomposable objects for
radical-square-zero linear Nakayama algebras, and exhibit an
infinite-order class modulo AR relations for the Kronecker algebra.

\section{Preliminaries}\label{sec:preliminaries}
In this section we recall the definitions and basic properties used in
the proofs, and fix our notation.

\subsection{Extriangulated categories and exact sequences}
\begin{definition}[{Nakaoka--Palu \cite[Definition 2.12]{NP}}]
An \emph{extriangulated category} is a triple
$(\mathcal B,\E,\mathfrak s)$ consisting of an additive category
$\mathcal B$, a biadditive functor
$\E:\mathcal B^{\mathrm{op}}\times\mathcal B\to\mathrm{Ab}$, and an
additive realization $\mathfrak s$ satisfying the axioms
\textup{(ET3)}, \textup{(ET3)$^{\mathrm{op}}$}, \textup{(ET4)} and
\textup{(ET4)$^{\mathrm{op}}$}. The realization assigns to each
$\delta\in\E(Z,X)$ an equivalence class of pairs of morphisms
$X\xrightarrow fY\xrightarrow gZ$. It respects direct sums and split
extensions and realizes morphisms of extensions. The axioms
\textup{(ET3)} and its dual complete compatible squares to morphisms
of these pairs, and \textup{(ET4)} and its dual give the corresponding
octahedral compatibility for their compositions.
\end{definition}

We call a realized extension an \emph{extriangle}, and write
$\delta:X\xrightarrow fY\xrightarrow gZ$ together with its extension
class $\delta\in\E(Z,X)$. The maps $f$ and $g$ are called an
\emph{inflation} and a \emph{deflation}, respectively. The extriangle
is split precisely when $\delta=0$, in which case its realization is
equivalent to $X\to X\oplus Z\to Z$. A full additive subcategory is
extension-closed if it contains the middle term whenever it contains
both the first and the third terms. The extension-closed subcategory then inherits the restricted extriangulated
structure by Nakaoka and Palu \cite[Remark 2.18]{NP}.

\begin{lemma}[{Nakaoka--Palu \cite[Corollary 3.12]{NP}}]\label{lem:HomE}
For every extriangle $\delta:X\xrightarrow fY\xrightarrow gZ$ in
$\mathcal B$ and every object $U\in\mathcal B$, the following sequences are exact:
\begin{equation}\label{eq:HomEcov}
\begin{aligned}
 \Hom(U,X)&\xrightarrow{f_*}\Hom(U,Y)
 \xrightarrow{g_*}\Hom(U,Z)\xrightarrow{\delta_\sharp}\E(U,X)\\
 &\xrightarrow{f_*}\E(U,Y)\xrightarrow{g_*}\E(U,Z),
\end{aligned}
\end{equation}
\begin{equation}\label{eq:HomEcontra}
\begin{aligned}
 \Hom(Z,U)&\xrightarrow{g^*}\Hom(Y,U)
 \xrightarrow{f^*}\Hom(X,U)\xrightarrow{\delta^\sharp}\E(Z,U)\\
 &\xrightarrow{g^*}\E(Y,U)\xrightarrow{f^*}\E(X,U).
\end{aligned}
\end{equation}
Here $\delta_\sharp(a)=a^*\delta$ and
$\delta^\sharp(b)=b_*\delta$.
\end{lemma}

\begin{definition}[{Nakaoka--Palu \cite[Definitions 3.23 and 3.25]{NP}}]
An object $P\in\mathcal B$ is \emph{projective} if every morphism $P\to Z$ lifts
through $g$ for every extriangle $X\to Y\xrightarrow gZ$.
An object $I\in\mathcal B$ is \emph{injective} if every morphism $X\to I$ extends
across $f$ for every extriangle $X\xrightarrow fY\to Z$.
The category has \emph{enough projectives} if every object $Z$
occurs in an extriangle $X\to P\to Z$ with $P$ projective; it has
\emph{enough injectives} if every object $X$ occurs in an extriangle
$X\to I\to Z$ with $I$ injective.
\end{definition}

\begin{lemma}[{Nakaoka--Palu \cite[Proposition 3.24]{NP}}]
\label{lem:proj-inj}
An object $P$ is projective if and only if $\E(P,V)=0$ for every
$V\in\mathcal B$. An object $I$ is injective if and only if
$\E(V,I)=0$ for every $V\in\mathcal B$.
\end{lemma}
Let $\mathcal B$ be an extriangulated category over $k$.
It is Hom-finite if $\dimk\Hom(U,V)<\infty$ for all $U,V\in\mathcal B$,
and $\E$-finite if $\dimk\E(U,V)<\infty$ for all $U,V\in\mathcal B$. In
the following, we will mainly concern with extriangulated categories which are Hom-finite, $\E$-finite and Krull-Schmidt.

\subsection{Extended module categories}
Let $A$ be a finite-dimensional algebra over an algebraic closed field $k$, and let $m\geq1$.
Let $\D=D^b(\modu A)$, with its standard $t$-structure, and let
$\D^{[a,b]}=\D^{\geq a}\cap\D^{\leq b}$. We write $\tau^{\leq q}$ and $\tau^{\geq q}$ for
standard cohomological truncation.

\begin{definition}[{\cite[Definition 2.2]{Zhou}}]
The $m$-extended module category of $A$ is
$\C_m(A)=\D^{[1-m,0]}$. Its morphisms are the morphisms of $\D$.
\end{definition}

For the rest of the article, we abbreviate $\C_m(A)$ to $\C$ when no
confusion can arise. Since $\C$ is closed under extensions and direct summands
in $\D$, it is an extriangulated category and the inherited extriangulated structure is
\[
 \E(Z,X)=\Hom_{\D}(Z,X[1]),
\]
with realization given by distinguished triangles whose first three terms
belong to $\C$, as described by Y.~Zhou \cite[Remark 1.4]{Zhou}.
The projective and injective objects are understood in the sense of
Lemma~\ref{lem:proj-inj} and we will use the following facts.

\begin{proposition}[{\cite[Lemma 3.1]{Zhou}}]\label{prop:structure}
The category $\C$ is $k$-linear, Hom-finite, $\E$-finite and Krull--Schmidt.
It has enough projectives and enough injectives, and
\[
 \operatorname{Proj}\C=\add A,\qquad
 \operatorname{Inj}\C=(\operatorname{inj}A)[m-1].
\]
\end{proposition}

We write
$\ind\C$ for a set of representatives of its nonzero indecomposable objects.

\begin{lemma}\label{lem:window}
For $U,V\in\C$ and $j\geq m$, we have $\Hom_{\D}(U,V[-j])=0$.
Moreover, for $0\leq i<m$,
\[
 \Hom_{\D}(A[i],U)\simeq\Hom_A(A,H^{-i}(U)),
\]
and, for an injective module $J$ and any integer $t$,
\[
 \Hom_{\D}(U,J[t])\simeq\Hom_A(H^{-t}(U),J).
\]
\end{lemma}
\begin{proof}
For the first assertion, $U\in\D^{\leq0}$ and
$V[-j]\in\D^{\geq1-m+j}\subseteq\D^{\geq1}$, so the $t$-structure
orthogonality applies. Since $A$ is projective, the functor
$\Hom_A(A,-)$ computes its derived Hom and commutes with cohomology.
This gives the second assertion. The last assertion follows in the same
way from exactness of $\Hom_A(-,J)$ and the fact that a stalk injective
complex computes derived Hom in the second variable.
\end{proof}

\subsection{Auslander--Reiten extriangles}
\begin{definition}[{Y.~Zhou \cite[Definition 3.3]{Zhou}}]
An extriangle
\begin{equation}\label{eq:triangle}
 \delta:\quad X\xrightarrow fY\xrightarrow gZ\xrightarrow\delta X[1]
\end{equation}
is an Auslander--Reiten extriangle (AR extriangle for short) if $\delta\ne0$, every morphism
$X\to X'$ in $\C$ which is not a section factors through $f$, and every
morphism $Z'\to Z$ in $\C$ which is not a retraction factors through $g$.
\end{definition}
The following proposition shows that in $\C$, AR extriangles always exist.

\begin{proposition}[{\cite[Theorem 3.12]{Zhou}}]\label{prop:ARexists}
Every non-projective indecomposable object of $\C$ is the endpoint of an
AR extriangle. Dually, every non-injective indecomposable
object is the starting point of an AR extriangle.
\end{proposition}

Iyama, Nakaoka and Palu \cite[Proposition 2.6]{INP} proved that
AR extriangles with a fixed endpoint are unique up to isomorphism. Their almost split properties are relative
to $\C$; an AR extriangle here need not be an AR triangle in $\D$.

\subsection{Grothendieck groups and local finiteness}
The split Grothendieck group of $\C$ is generated by object classes modulo
$[X\oplus Y]=[X]+[Y]$. Since $\C$ is Krull--Schmidt,
\[
 \ksp(\C)=\bigoplus_{U\in\ind\C}\mathbb Z[U].
\]
For an extriangle
\[
 \delta:\quad X\xrightarrow fY\xrightarrow gZ
\]
in $\C$, define its relation by $r_\delta=[X]-[Y]+[Z]$ in
$\ksp(\C)$. Set
\[
 \Ex(\C)=\langle r_\delta:\delta\text{ an extriangle}\rangle_{\mathbb Z},
 \qquad
 \AR(\C)=\langle r_\delta:\delta\text{ an AR extriangle}\rangle_{\mathbb Z}.
\]
Then $\AR(\C)\subseteq\Ex(\C)$ and
$K_0(\C)=\ksp(\C)/\Ex(\C)$; see \cite[Section 4]{ZZ}.
For a module category, the relation subgroup of its split
Grothendieck group is generated by short exact sequences; for a
triangulated category, it is generated by distinguished triangles.

We call $\C$ \emph{finite}, or of \emph{finite representation type},
if $\ind\C$ is finite. Following Zhu and Zhuang
\cite[Definition 3.1]{ZZ}, $\C$ is locally finite if for each
$V\in\ind\C$ both sets
\[
 \{U\in\ind\C:\Hom(U,V)\ne0\},\qquad
 \{U\in\ind\C:\Hom(V,U)\ne0\}
\]
are finite.

The same definitions of $\ksp$, $\Ex$, $\AR$ and local finiteness
will be used for any extriangulated category $\mathcal B$.

\begin{theorem}[{Zhu--Zhuang \cite[Theorem 4.4]{ZZ}}]\label{thm:ZZ}
Let $\mathcal B$ be a skeletally small, $k$-linear, Hom-finite,
$\E$-finite Krull--Schmidt extriangulated category. If $\mathcal B$ is
locally finite, then $\Ex(\mathcal B)=\AR(\mathcal B)$.
\end{theorem}

\section{A criterion with bounded negative extensions}\label{sec:general}
Let $\mathcal B$ be a skeletally small, $k$-linear, Hom-finite,
$\E$-finite Krull--Schmidt extriangulated category. For an extriangle
$\delta:X\to Y\xrightarrow g Z$, its right defect is
\[
 D_\delta(U)=\coker\bigl(\Hom(U,Y)\xrightarrow{g_*}\Hom(U,Z)\bigr).
\]
An \emph{extension generator} is an object $F\in\mathcal B$ such that every object
is obtained from $\add F$ by finitely many extensions and direct
summands. The number of extensions may depend on the object.

\begin{definition}\label{def:negative}
A bounded negative extension sequence for the second variable consists
of additive bifunctors
$N^{-j}:\mathcal B^{\mathrm{op}}\times\mathcal B\to\modu k$, $j\geq0$,
with $N^0=\Hom_{\mathcal B}$, such that $N^{-j}=0$ for $j\geq b$, for some $b\geq1$,
and every extriangle $\delta:X\to Y\to Z$ induces the exact sequence
\begin{equation}\label{eq:negativeexact}
\begin{aligned}
0&\longrightarrow N^{1-b}(U,X)\longrightarrow N^{1-b}(U,Y)
 \longrightarrow N^{1-b}(U,Z)\longrightarrow\cdots\\
 &\longrightarrow\Hom(U,X)\longrightarrow\Hom(U,Y)
 \longrightarrow\Hom(U,Z)\longrightarrow D_\delta(U)\longrightarrow0
\end{aligned}
\end{equation}
for every $U\in\mathcal B$.
\end{definition}

The definition above specifies a finite covariant $\delta$-functor
in the second variable in the framework of Gorsky, Nakaoka and Palu
\cite[Section 4]{GNP}, with a uniform bound $b$ and finite-dimensional
values. Adachi, Enomoto and Tsukamoto \cite[Definition 2.3]{AET} introduced a negative first extension
structure supplying degree $-1$ to degree $0$ exactness in both
variables. Here we use only the second-variable sequence and require
exactness through all negative degrees down to the stated zero
boundary. Ambient negative Hom in $D^b(\modu A)$ provides these data
for $\C_m(A)$. In an exact category one can take $b=1$ and all
negative bifunctors are zero.

Now we define a homomorphism on the split Grothendieck group by
\begin{equation}\label{eq:chi}
 \chi^-_U(V)=\sum_{j=0}^{b-1}(-1)^j\dimk N^{-j}(U,V),\qquad
 \chi^-_U:\ksp(\mathcal B)\longrightarrow\mathbb Z.
\end{equation}

\begin{lemma}\label{lem:defect}
For every extriangle $\delta$ and object $U$,
$\chi^-_U(r_\delta)=\dimk D_\delta(U)$.
\end{lemma}
\begin{proof}
Take alternating dimensions in \eqref{eq:negativeexact}, starting at
its right endpoint. The three terms in degree $-j$ contribute
$(-1)^j(\dimk N^{-j}(U,X)-\dimk N^{-j}(U,Y)+\dimk N^{-j}(U,Z))$.
Summing gives the assertion, since the far-left term preceding
$N^{1-b}(U,X)$ is zero.
\end{proof}

Let
\[
 \operatorname{Ker}_\delta(U)=\ker(\Hom(U,X)\to\Hom(U,Y)).
\]
By Lemma~\ref{lem:HomE}, the ordinary three-term Hom-dimension
difference is explicitly
\[
\begin{aligned}
 &\dimk\Hom(U,X)-\dimk\Hom(U,Y)+\dimk\Hom(U,Z)\\
 &\hspace{25mm}=\dimk D_\delta(U)+\dimk\operatorname{Ker}_\delta(U).
\end{aligned}
\]
The sign convention agrees with $r_\delta=[X]-[Y]+[Z]$.
Taking alternating dimensions in the negative part of
\eqref{eq:negativeexact} gives
\[
\begin{aligned}
 &\sum_{j=1}^{b-1}(-1)^j\bigl(\dimk N^{-j}(U,X)
       -\dimk N^{-j}(U,Y)+\dimk N^{-j}(U,Z)\bigr)\\
 &\hspace{25mm}=-\dimk\operatorname{Ker}_\delta(U).
\end{aligned}
\]
Thus the negative terms cancel the additional left boundary term.
This is the term that allows a second support point in the ordinary
Hom calculation of Haugland \cite[Lemma 2.3]{Haugland}. In the following of this section, we assume that $\mathcal B$ has AR extriangles.

\begin{lemma}\label{lem:ARsupport}
Let $\delta_Z$ be an AR extriangle ending at $Z\in\mathcal B$. Set
\[
 d_Z=\dimk(\End(Z)/\rad\End(Z)).
\]
For indecomposable $U$,
\[
 \chi^-_U(r_{\delta_Z})=
 \begin{cases}d_Z,&U\simeq Z,\\0,&U\not\simeq Z.\end{cases}
\]
\end{lemma}
\begin{proof}
For $U\not\simeq Z$, every map $U\to Z$ is a non-retraction and
factors through the right almost split map, so $D_{\delta_Z}(U)=0$.
For $U=Z$, its image in $\End(Z)$ consists precisely of the nonunits:
all nonunits factor, and factoring a unit would split the map.
Hence $D_{\delta_Z}(Z)=\End(Z)/\rad\End(Z)$. Apply
Lemma~\ref{lem:defect}.
\end{proof}

\begin{theorem}\label{thm:general}
Assume that $\mathcal B$ has a bounded negative extension sequence,
enough injectives, an extension generator $F$, finitely many
indecomposable projectives, and AR extriangles. Choose an extriangle
$\epsilon:F\to I\to T$ with $I$ injective. Then the following
conditions are equivalent:
\begin{enumerate}
\item $\mathcal B$ is finite;
\item $\mathcal B$ is locally finite;
\item $\Ex(\mathcal B)=\AR(\mathcal B)$;
\item $r_\epsilon\in\AR(\mathcal B)$.
\end{enumerate}
\end{theorem}
\begin{proof}
We first identify the value of $\chi^-_U$ on the chosen relation.
Applying Lemma~\ref{lem:HomE} to $\epsilon$ gives the exact sequence
\[
\begin{aligned}
 \Hom(U,F)&\longrightarrow\Hom(U,I)\longrightarrow\Hom(U,T)
 \xrightarrow{\epsilon_\sharp}\E(U,F)\\
 &\longrightarrow\E(U,I)\longrightarrow\E(U,T).
\end{aligned}
\]
Since $I$ is injective, Lemma~\ref{lem:proj-inj} gives $\E(U,I)=0$.
Consequently $\epsilon_\sharp$ induces an isomorphism
\[
 D_\epsilon(U)=\coker(\Hom(U,I)\to\Hom(U,T))\simeq\E(U,F).
\]
By Lemma~\ref{lem:defect},
\begin{equation}\label{eq:generatorvalue}
 \chi^-_U(r_\epsilon)=\dimk\E(U,F).
\end{equation}

For fixed $U$, the full subcategory
$\mathcal V_U=\{V:\E(U,V)=0\}$ is closed under finite direct sums
, direct summands and extension. If $\E(U,F)=0$,
$\mathcal V_U$ contains $\add F$, and hence every object by the
extension-generator hypothesis. Lemma~\ref{lem:proj-inj} then makes
$U$ projective. Conversely, projectivity immediately gives
$\E(U,F)=0$. Thus \eqref{eq:generatorvalue} is strictly positive
for every non-projective indecomposable $U$.

\emph{(4) implies (1).} Membership in the subgroup $\AR(\mathcal B)$
means that there are finitely many AR extriangles
$\delta_1,\ldots,\delta_s$, with endpoints $Z_1,\ldots,Z_s$, and
integers $a_1,\ldots,a_s$ such that
\[
 r_\epsilon=\sum_{i=1}^s a_i r_{\delta_i}.
\]
If an indecomposable $U$ is not isomorphic to any $Z_i$, applying
$\chi^-_U$ and Lemma~\ref{lem:ARsupport} yields
\[
 \dimk\E(U,F)=\chi^-_U(r_\epsilon)
       =\sum_{i=1}^s a_i\chi^-_U(r_{\delta_i})=0.
\]
Such an object must be projective by the preceding paragraph.
Therefore every non-projective indecomposable is isomorphic to one
of $Z_1,\ldots,Z_s$. There are finitely many projective
indecomposables by hypothesis, so $\mathcal B$ is finite. The same
argument covers $s=0$, when all objects are projective.

\emph{(1) and (2) are equivalent.} Finiteness plainly implies local
finiteness because each Hom support is a subset of $\ind\mathcal B$.
For the converse, suppose first that $\Hom(F,U)=0$. The class
$\mathcal W_U=\{V:\Hom(V,U)=0\}$ contains $\add F$ and is closed
under dircet sum, dircet summands and extensions.
In particular $\End(U)=0$, whence $1_U=0$ and $U$ is a zero object.
It follows that $\Hom(F,U)\ne0$ for every nonzero $U$.
Write $F=\bigoplus_{a=1}^t F_a$ as a finite sum of indecomposables,
allowing repetitions. Then
\[
 \ind\mathcal B\subseteq
 \bigcup_{a=1}^t\{U\in\ind\mathcal B:\Hom(F_a,U)\ne0\}.
\]
Under local finiteness each set on the right is finite, proving (1).
If $F=0$, the same argument shows that $\mathcal B$ is the zero category.

Finally, (2) implies (3) by Theorem~\ref{thm:ZZ}, and (3) implies
(4) because $r_\epsilon$ is an extriangle relation and hence belongs
to $\Ex(\mathcal B)$. These implications prove all four equivalences.
\end{proof}

\begin{corollary}\label{cor:ARbasis}
Under the hypotheses of Theorem~\ref{thm:general}, choose one AR
extriangle for each endpoint for which such an
extriangle exists. Their relations are linearly independent. If
$\mathcal B$ is of finite type, they are indexed by all
non-projective indecomposable objects and form a $\mathbb Z$-basis
of $\Ex(\mathcal B)$. Every extriangle then has the unique expansion
\begin{equation}\label{eq:ARcoefficients}
 r_\delta=\sum_{Z\text{ non-projective}}
       \frac{\dimk D_\delta(Z)}{d_Z}\,r_{\delta_Z},
\end{equation}
with nonnegative integral coefficients.
\end{corollary}
\begin{proof}
Let the endpoints pairwise nonisomorphic. Suppose that a finite
integral linear combination of their relations is zero:
\[
 \sum_{i=1}^s a_i r_{\delta_{Z_i}}=0.
\]
For each $j$, evaluation by $\chi^-_{Z_j}$ and
Lemma~\ref{lem:ARsupport} gives $a_jd_{Z_j}=0$. Since the local
endomorphism ring of $Z_j$ has a nonzero division-ring quotient,
$d_{Z_j}$ is a positive integer, so $a_j=0$. This proves independence,
including when the full family of endpoints is infinite.

Now assume that $\mathcal B$ is finite. It is locally finite, so
Zhu and Zhuang \cite[Theorem 3.12]{ZZ} showed that $\mathcal B$ has AR extriangles. By Iyama, Nakaoka and Palu
\cite[Proposition 2.6]{INP}, two AR extriangles with the same endpoint
are isomorphic and give the same relation. Theorem~\ref{thm:general}
gives $\Ex(\mathcal B)=\AR(\mathcal B)$, so the chosen independent
relations span $\Ex(\mathcal B)$ over $\mathbb Z$.

In particular, every extriangle $\delta$ has a unique expression
$r_\delta=\sum_Z a_Zr_{\delta_Z}$ with $a_Z\in\mathbb Z$.
Applying Lemmas~\ref{lem:defect} and \ref{lem:ARsupport} at a fixed
non-projective indecomposable $W$ gives
\[
 \dimk D_\delta(W)=\chi^-_W(r_\delta)
    =\sum_Z a_Z\chi^-_W(r_{\delta_Z})=a_Wd_W.
\]
Therefore $a_W=\dimk D_\delta(W)/d_W$, as asserted in
\eqref{eq:ARcoefficients}. Integrality follows from the integral
spanning statement just proved; nonnegativity follows from
$\dimk D_\delta(W)\geq0$ and $d_W>0$.
\end{proof}

Formula \eqref{eq:ARcoefficients} has the same basis-expansion form as
\cite[Theorem 4.13 and Corollary 4.14]{PPPP}; the functional here is
the right defect evaluated through the bounded negative sequence.

\begin{remark}\label{rem:hearts}
Let $\mathcal T$ be a triangulated category with a bounded t-structure whose heart is a length category with finitely many simple objects. The sum of their shifts $0,\ldots,m-1$ extension-generates
$\mathcal T^{[1-m,0]}$, and ambient negative Hom vanishes in degrees
$\leq-m$. Thus the criterion also applies when this extended heart
satisfies the remaining finiteness and injectivity hypotheses of
Theorem~\ref{thm:general}. For broader constructions and tilting
theory of extended hearts, see \cite{MNO,WZ}.
\end{remark}

\section{Grothendieck groups in Extended module categories}\label{sec:module}
We now apply the general criterion in Section~\ref{sec:general} to $\C=\C_m(A)$. We first record
the cohomological description of $K_0(\C)$ and its relation subgroup.

\subsection{The cohomological Euler characteristic}
\begin{proposition}\label{thm:K0}
The inclusions $\modu A\hookrightarrow\C_m(A)\hookrightarrow\D$ induce
isomorphisms
\[
 K_0(\modu A)\xrightarrow{j_*}K_0(\C_m(A))
 \xrightarrow{\ell_*}K_0(\D).
\]
The map
\begin{equation}\label{eq:eu}
 \eu_m:K_0(\C_m(A))\longrightarrow K_0(\modu A),\qquad
 [X]\longmapsto\sum_{i=0}^{m-1}(-1)^i[H^{-i}(X)]
\end{equation}
is well-defined and is inverse to $j_*$.
\end{proposition}
\begin{proof}
Every short exact sequence of modules gives an extriangle in the
degree-zero heart, so $j_*$ is well-defined. The inclusion into $\D$
preserves the defining relations as well, and therefore induces $\ell_*$.

We next show that \eqref{eq:eu} respects relations. For an extriangle
$X\to Y\to Z$ in $\C_m(A)$, the cohomology long exact sequence is finite:
\[
\begin{aligned}
0&\longrightarrow H^{1-m}(X)\longrightarrow H^{1-m}(Y)
 \longrightarrow H^{1-m}(Z)\longrightarrow\cdots\\
 &\longrightarrow H^0(X)\longrightarrow H^0(Y)
 \longrightarrow H^0(Z)\longrightarrow0.
\end{aligned}
\]
For any finite exact sequence of modules, its alternating sum
vanishes in $K_0(\modu A)$: this follows by breaking it into short exact
sequences using the images of its maps, whose classes cancel in adjacent
positions. Applying this to the displayed sequence gives
\[
 \sum_{q=1-m}^0(-1)^q
       \bigl([H^q(X)]-[H^q(Y)]+[H^q(Z)]\bigr)=0.
\]
Thus $\eu_m$ kills every defining relation and is well-defined.

For a module $M$ in degree zero, $\eu_mj_*([M])=[M]$.
To prove the other composition is the identity, let $X\in\C$.
For $1-m\leq q\leq0$, the standard truncation triangles are
\begin{equation}\label{eq:trunc}
 \tau^{\leq q-1}X\longrightarrow\tau^{\leq q}X
 \longrightarrow H^q(X)[-q]\longrightarrow\tau^{\leq q-1}X[1].
\end{equation}
Their first three terms lie in $\C$, including
$\tau^{\leq-m}X=0$. Successively applying their relations gives
\begin{equation}\label{eq:filtrationclass}
 [X]=\sum_{i=0}^{m-1}[H^{-i}(X)[i]]\quad\text{in }K_0(\C).
\end{equation}
For any module $M$ and $0\leq i<m-1$, the triangle
\[
 M[i]\longrightarrow0\longrightarrow M[i+1]
 \xrightarrow{1}M[i+1]
\]
also has its first three terms in $\C$. Hence
$[M[i+1]]=-[M[i]]$, and induction gives
\begin{equation}\label{eq:sign}
 [M[i]]=(-1)^i[M]\qquad(0\leq i<m).
\end{equation}
Equations \eqref{eq:filtrationclass} and \eqref{eq:sign} show
$j_*\eu_m([X])=[X]$. For $m=1$, the shift step is empty and the same
argument applies. This proves the first isomorphism.

For completeness, the same finite cohomology argument defines a map on
$K_0(\D)$ by
\[
 \eu_{\D}([X])=\sum_{q\in\mathbb Z}(-1)^q[H^q(X)].
\]
Each object is bounded, so this sum is finite. Truncation triangles and
$[V[1]]=-[V]$ show that $\eu_{\D}$ is inverse to the homomorphism induced
by $\modu A\hookrightarrow\D$. This latter homomorphism is
$\ell_*j_*$. Since $j_*$ and $\ell_*j_*$ are isomorphisms, so is
$\ell_*$. Moreover $\eu_{\D}\ell_*=\eu_m$.
\end{proof}

\begin{corollary}\label{cor:simples}
The classes of the simple right $A$-modules $S_1,\ldots,S_n$
form a $\mathbb Z$-basis of $K_0(\C_m(A))$.
For $m\leq r$, the inclusion $\C_m(A)\hookrightarrow\C_r(A)$ induces
an isomorphism, which is the identity under the above identifications
with $K_0(\modu A)$.
\end{corollary}
\begin{proof}
Every finitely generated $A$-module $M$ has a finite composition
series. Its successive short exact sequences give
\[
 [M]=\sum_{a=1}^n[M:S_a][S_a]\quad\text{in }K_0(\modu A),
\]
where $[M:S_a]$ is the number of composition factors isomorphic to
$S_a$. This number is independent of the series by the
Jordan--H\"older theorem. Thus the simple classes generate $K_0(\modu A)$.
The functions $[M:S_a]$ are additive on short exact sequences, so
they define a homomorphism
\[
 K_0(\modu A)\longrightarrow\mathbb Z^n,\qquad
 [M]\longmapsto([M:S_1],\ldots,[M:S_n]).
\]
It sends $[S_a]$ to the $a$-th standard basis vector. Applying it to
any integral relation among the $[S_a]$ forces all coefficients to
vanish. The first assertion now follows from the isomorphism $j_*$
in Proposition~\ref{thm:K0}.

For $m\leq r$, let $\iota:\C_m(A)\hookrightarrow\C_r(A)$ be the
inclusion. It preserves extriangles because both structures come
from the same derived category. It therefore induces $\iota_*$ on
Grothendieck groups. For $X\in\C_m(A)$ we have $H^{-i}(X)=0$ for
$i\geq m$, and hence
\[
 \eu_r\iota_*([X])=\sum_{i=0}^{r-1}(-1)^i[H^{-i}(X)]
                 =\eu_m([X]).
\]
Since both Euler maps are isomorphisms by Proposition~\ref{thm:K0},
$\iota_*=\eu_r^{-1}\eu_m$ is an isomorphism and becomes the identity
after identifying both groups with $K_0(\modu A)$.
\end{proof}

The preceding proof also determines the relation subgroup without any
finiteness assumption on $\ind\C$. For $U\in\ind\C$, define
its Euler coordinates by
\[
 e_a(U)=\sum_{i=0}^{m-1}(-1)^i[H^{-i}(U):S_a],\qquad
 e(U)=(e_1(U),\ldots,e_n(U)).
\]
Here $[H^{-i}(U):S_a]$ is the composition multiplicity of $S_a$ in
the finite-length module $H^{-i}(U)$, not a Hom dimension. Thus
$e_a(U)$ is the alternating sum of those multiplicities across the
cohomological degrees of $U$. Equivalently, $e(U)$ is the coordinate
vector of $\eu_m([U])$ in the simple basis of $K_0(\modu A)$.

\begin{proposition}\label{prop:kernel}
Let $\mathscrS=\{S_1,\ldots,S_n\}\subseteq\ind\C$. For
$U\in\ind\C\setminus\mathscrS$, set
\[
 \rho_U=[U]-\sum_{a=1}^n e_a(U)[S_a].
\]
Then
\begin{equation}\label{eq:kernelbasis}
 \Ex(\C)=\bigoplus_{U\in\ind\C\setminus\mathscrS}\mathbb Z\rho_U.
\end{equation}
In particular, there is a split exact sequence of abelian groups
\[
0\longrightarrow\Ex(\C)\longrightarrow\ksp(\C)
 \xrightarrow{\ e\ }\mathbb Z^n\longrightarrow0.
\]
If $N=\#\ind\C$ is finite, then $\operatorname{rank}\Ex(\C)=N-n$.
\end{proposition}
\begin{proof}
Let $\pi:\ksp(\C)\to K_0(\C)$ be the quotient map. Compose it
with $\eu_m$ from Proposition~\ref{thm:K0} and the simple-coordinate
isomorphism from Corollary~\ref{cor:simples}. The resulting map is
exactly $e$, because its value on $[U]$ is $(e_1(U),\ldots,e_n(U))$.
Consequently $\ker e=\ker\pi=\Ex(\C)$.

Since $e([S_a])$ is the $a$-th standard basis vector, $e(\rho_U)=0$.
Thus each $\rho_U$ belongs to $\Ex(\C)$. To prove that they span,
write an arbitrary $x\in\ker e$ as the finite sum
\[
 x=\sum_{U\notin\mathscrS}c_U[U]+\sum_{a=1}^n b_a[S_a].
\]
Substituting $[U]=\rho_U+\sum_a e_a(U)[S_a]$ gives
\[
 x=\sum_{U\notin\mathscrS}c_U\rho_U+
       \sum_{a=1}^n\left(b_a+\sum_{U\notin\mathscrS}c_Ue_a(U)\right)[S_a].
\]
The parenthesized coefficients are precisely the coordinates of
$e(x)$, so they vanish. This proves spanning. If
$\sum_U c_U\rho_U=0$, comparing the coefficient of the non-simple
indecomposable class $[U]$ in the original free basis of $\ksp(\C)$
gives $c_U=0$ for every $U$. Hence the spanning family is independent,
and \eqref{eq:kernelbasis} follows. All sums involved in an individual
element are finite, so this proof also covers infinitely many
indecomposables.

The homomorphism $\sigma:\mathbb Z^n\to\ksp(\C)$ sending the
$a$-th standard basis vector to $[S_a]$ satisfies $e\sigma=1$.
It splits the asserted exact sequence. If $\ind\C$ has $N$ elements,
the basis in \eqref{eq:kernelbasis} has exactly $N-n$ elements,
which gives the stated rank.
\end{proof}

\begin{remark}
The rank of $K_0(\C)$ is independent of $m$; the rank of its relation
subgroup can be infinite.
\end{remark}

\subsection{The explicit relation and the main criterion}
Take $N^{-j}(U,V)=\Hom_{\D}(U,V[-j])$ and $b=m$. The ambient
long exact Hom sequence, together with Lemma~\ref{lem:window}, gives
\eqref{eq:negativeexact}. Thus the function of \eqref{eq:chi} is
\[
 \chi^-_U(V)=\sum_{j=0}^{m-1}(-1)^j\dimk\Hom_{\D}(U,V[-j]).
\]
Let $S=\bigoplus_aS_a$ and $F=\bigoplus_{i=0}^{m-1}S[i]$.

\begin{lemma}\label{lem:detector}
The object $F$ is an extension generator. The extriangle
$\epsilon:F\to J[m-1]\to T$ in \eqref{eq:introdetector} has injective
middle term, and
\[
 \chi^-_U(r_\epsilon)=\dimk\E(U,F),\qquad
 \E(U,F)=0\ \Longleftrightarrow\ U\in\add A.
\]
\end{lemma}
\begin{proof}
For a module $M$, each short exact sequence in a composition series,
shifted by $i$ with $0\leq i<m$, gives an extriangle in $\C$.
Thus $M[i]$ is obtained by finitely many extensions of the objects
$S_a[i]$, all of which belong to $\add F$. For a general object
$X\in\C$, the truncation triangles \eqref{eq:trunc} successively
build $X$ from the shifted cohomology modules $H^{-i}(X)[i]$.
Combining these two finite constructions proves that $F$ is an
extension generator.

We describe the chosen extriangle explicitly. For each
$0\leq i<m-1$, take the distinguished triangle
\[
 S[i]\longrightarrow0\longrightarrow S[i+1]
       \xrightarrow{1}S[i+1].
\]
Also shift the triangle associated to
$0\to S\to J\to Q\to0$ by $m-1$:
\[
 S[m-1]\longrightarrow J[m-1]\longrightarrow Q[m-1]
           \longrightarrow S[m].
\]
Their direct sum has first three terms $F,I,T$, where
\begin{equation}\label{eq:fixedtriangle}
 I=J[m-1],\qquad T=\bigoplus_{i=1}^{m-1}S[i]\oplus Q[m-1].
\end{equation}
These three terms belong to $\C$, so the direct-sum triangle induces
the extriangle $\epsilon$ in \eqref{eq:introdetector}. When $m=1$
there are no shift triangles, the sum over $1\leq i\leq m-1$ is
the zero object, and the construction is simply $S\to J\to Q$.
Proposition~\ref{prop:structure} identifies $I$ as injective.

Lemma~\ref{lem:HomE} now gives the exact sequence
\[
\begin{aligned}
 \Hom(U,F)&\longrightarrow\Hom(U,J[m-1])\longrightarrow\Hom(U,T)\\
 &\longrightarrow\E(U,F)\longrightarrow\E(U,J[m-1])=0.
\end{aligned}
\]
The zero follows from Lemma~\ref{lem:proj-inj}. It identifies
$D_\epsilon(U)$ with $\E(U,F)$. Lemma~\ref{lem:defect} therefore gives
\[
 \chi^-_U(r_\epsilon)=\dimk D_\epsilon(U)=\dimk\E(U,F).
\]
Finally, if $\E(U,F)=0$, additivity and the exact segment
$\E(U,V_1)\to\E(U,V)\to\E(U,V_2)$ in Lemma~\ref{lem:HomE}
propagate this vanishing from $\add F$ to all objects built by
extensions and summands. Since $F$ is an extensioin generator, $\E(U,-)=0$.
Lemma~\ref{lem:proj-inj} and Proposition~\ref{prop:structure} imply
$U\in\add A$. Conversely, $U\in\add A$ is projective and so
$\E(U,F)=0$.
\end{proof}

\begin{theorem}\label{thm:finite}
For $\epsilon$ in \eqref{eq:introdetector}, the following are equivalent:
\begin{enumerate}
\item $\C_m(A)$ has finitely many indecomposable objects up to isomorphism;
\item $\C_m(A)$ is locally finite;
\item $\Ex(\C_m(A))=\AR(\C_m(A))$;
\item $r_\epsilon\in\AR(\C_m(A))$.
\end{enumerate}
\end{theorem}
\begin{proof}
Proposition~\ref{prop:structure} gives Hom and $\E$ finiteness,
the Krull--Schmidt property, enough injectives, and finitely many
indecomposable projectives. The negative sequence above and
Lemma~\ref{lem:detector} verify the other hypotheses of
Theorem~\ref{thm:general}, which proves the assertions.
\end{proof}

\begin{corollary}\label{cor:modulebasis}
One AR relation for each non-projective indecomposable gives a linearly
independent family in $\ksp(\C_m(A))$. If $N=\#\ind\C_m(A)$ is finite,
this is a $\mathbb Z$-basis of $\Ex(\C_m(A))$ with $N-n$ elements, where $n$
is the number of simple modules, and \eqref{eq:ARcoefficients} gives
the coefficients of every extriangle relation.
\end{corollary}
\begin{proof}
Proposition~\ref{prop:ARexists} supplies an AR extriangle at every
non-projective indecomposable endpoint. The bounded negative Hom
sequence established above allows us to apply
Corollary~\ref{cor:ARbasis}, which proves linear independence.
If $N$ is finite, Theorem~\ref{thm:finite} gives $\Ex=\AR$, so
these relations form an integral basis and their coefficients are
given by \eqref{eq:ARcoefficients}. By
Proposition~\ref{prop:structure}, the projective objects are $\add A$.
There are $n$ indecomposable projective modules up to isomorphism,
one projective cover for each simple module. Thus the basis has
$N-n$ elements, in agreement with Proposition~\ref{prop:kernel}.
\end{proof}

\begin{corollary}\label{cor:hereditary}
If $A$ is hereditary, the conditions of Theorem~\ref{thm:finite} are
equivalent to $A$ being representation-finite. If $A$ is hereditary and
representation-finite, with $t$ indecomposable modules and $n$ simple
modules, then
\[
 \#\ind\C_m(A)=mt,\qquad
 \operatorname{rank}\Ex(\C_m(A))=mt-n.
\]
\end{corollary}
\begin{proof}
We recall why the indecomposables are the objects $M[i]$, where
$M\in\ind(\modu A)$ and $0\leq i<m$; see also
\cite[Example 3.15]{Zhou}. Inductively split the truncation triangles
\eqref{eq:trunc}. After the preceding truncation has been decomposed,
the components of its connecting map lie in
\[
 \Hom_{\D}(H^q(X)[-q],H^r(X)[1-r])
 =\Ext_A^{q-r+1}(H^q(X),H^r(X)),\qquad r<q.
\]
These groups vanish since their degrees are at least two. Hence
$X\simeq\bigoplus_iH^{-i}(X)[i]$.
Different shifts of nonzero modules are distinguished by cohomology,
and indecomposability is preserved by shifting. The assertions follow
from Theorem~\ref{thm:finite} and Proposition~\ref{prop:kernel}.
\end{proof}

\begin{remark}\label{rem:notconverse}
For every finite-dimensional algebra $A$, if $\C_m(A)$ is of finite
type, then $A$ is representation-finite, by the fully faithful inclusion of
$\modu A$. The converse fails: Rundsveen
\cite[Lemma 4.11(4), Theorem 4.1]{Rundsveen} proved that
$\C_2(\Lambda(9,6))$ is infinite, where $\Lambda(n,\ell)$ is the linear
Nakayama algebra on $n$ vertices with all paths of length $\ell$ zero. The algebra itself is
representation-finite. His left-module convention is transferred here by
taking the opposite algebra and reversing the vertex order.
\end{remark}

\subsection{Two-term silting categories}
We next recover the two-term finiteness criterion of PPPP
\cite{PPPP} for
$K^{[-1,0]}(\operatorname{proj}A)$. Their more general presentation
categories are built from triangles $P_1\to P_0\to X\to P_1[1]$
with $P_0,P_1\in\add P$ for a rigid object $P$; see PPPP
\cite[Section 4.7.2]{PPPP}. Their proof of the finiteness criterion
refers to PPPP \cite[proof of Proposition 3.22]{PPPP}, where the
relation $[P]+[P[1]]$ is already used. We explain how the criterion
of Section~\ref{sec:general} recovers the same four equivalent
conditions under a bounded negative-Hom hypothesis.

\begin{definition}[{Aihara--Iyama \cite[Section 2]{AI}}]\label{def:silting}
An object $P$ in a triangulated category $\mathcal T$ is
\emph{silting} if
\[
 \Hom_{\mathcal T}(P,P[t])=0\quad(t>0),\qquad
 \operatorname{thick}(P)=\mathcal T.
\]
Here $\operatorname{thick}(P)$ is the smallest full triangulated
subcategory containing $P$ and closed under direct summands.
\end{definition}

For the remainder of this subsection, let $\mathcal T$ be an
idempotent-complete, skeletally small, $k$-linear triangulated
category and let $P$ be a silting object. Assume that
$\Hom_{\mathcal T}(P,P[t])$ is finite-dimensional for every integer
$t$, and that $\Hom_{\mathcal T}(P,P[-j])=0$ for all sufficiently
large $j$. Set
\[
 \mathcal B=\add P*\add P[1]
   =\{X:P_0\to X\to P_1[1]\to P_0[1],\ P_0,P_1\in\add P\}.
\]
This category has the following inherited structure.

\begin{lemma}\label{lem:silting-structure}
The category $\mathcal B$ is extension-closed and closed under
direct summands in $\mathcal T$. It is a Hom-finite, $\E$-finite
Krull--Schmidt extriangulated category with
$\E(Z,X)=\Hom_{\mathcal T}(Z,X[1])$. Moreover,
\[
 \operatorname{Proj}\mathcal B=\add P,\qquad
 \operatorname{Inj}\mathcal B=\add P[1],
\]
and $\mathcal B$ has enough projectives and enough injectives.
\end{lemma}
\begin{proof}
The subcategory of objects $Y$ for which $\Hom(P,Y[t])$ is finite
for every $t$ is thick, by the long exact Hom sequence and closure
under summands. It contains $P$, so it is all of $\mathcal T$.
Repeating this argument in the first variable shows that
$\Hom(X,Y[t])$ is finite for all $X,Y,t$. Idempotent completeness
then gives the Krull--Schmidt property.

Put $\mathcal P=\add P$. Silting orthogonality implies
$\mathcal P[1]*\mathcal P=\add(P\oplus P[1])$: the connecting map
of such a triangle lies in $\Hom(\mathcal P,\mathcal P[2])=0$.
It also gives $\mathcal P*\mathcal P=\mathcal P$ and the shifted
analogue. Associativity of extension products, which follows from
the octahedral axiom, therefore gives
\[
 (\mathcal P*\mathcal P[1])*(\mathcal P*\mathcal P[1])
       \subseteq\mathcal P*\mathcal P[1].
\]
Hence $\mathcal B$ is extension-closed.

For summand closure, take $X=X'\oplus X''\in\mathcal B$ and a
defining triangle $P_0\xrightarrow aX\to P_1[1]\to P_0[1]$.
Applying $\Hom(P,-)$ shows that $a$ is a right $\add P$-approximation.
Hom-finiteness and the Krull--Schmidt property give minimal right
$\add P$-approximations of $X'$, $X''$, and $X$. The direct sum
of the first two is a minimal approximation of $X$, and any
approximation of $X$ is isomorphic to this one plus a zero map
from an additional summand of its source. Consequently the cones
of the minimal approximations of $X'$ and $X''$ are summands of
$P_1[1]$. They belong to $\add P[1]$, and their triangles show
$X',X''\in\mathcal B$. The inherited extriangulated structure now
follows from the extension-closed-subcategory construction recalled
in Section~\ref{sec:preliminaries}. Hom and $\E$ finiteness follow
from the corresponding shifted Hom spaces in $\mathcal T$.

For $X\in\mathcal B$, its defining triangle and silting
orthogonality give $\Hom(P,X[1])=0$ and
$\Hom(X,P[2])=0$. Lemma~\ref{lem:proj-inj} therefore makes
$\add P$ projective and $\add P[1]$ injective in $\mathcal B$.
Rotating the defining triangle produces extriangles
\[
 P_1\longrightarrow P_0\longrightarrow X,
 \qquad X\longrightarrow P_1[1]\longrightarrow P_0[1].
\]
They give enough projectives and injectives. If $X$ is projective,
the first extriangle splits by Lemma~\ref{lem:proj-inj}, making
$X$ a summand of $P_0$. If $X$ is injective, the second splits,
making $X$ a summand of $P_1[1]$. This proves both equalities.
\end{proof}

\begin{corollary}[Two-term silting subcategories]\label{cor:silting}
Under the assumptions of this subsection, the following four
conditions are equivalent:
\begin{enumerate}
\item $\mathcal B$ is finite;
\item $\mathcal B$ is locally finite;
\item $\Ex(\mathcal B)=\AR(\mathcal B)$;
\item $[P]+[P[1]]\in\AR(\mathcal B)$.
\end{enumerate}
\end{corollary}
\begin{proof}
We verify the hypotheses of Theorem~\ref{thm:general} and identify
its distinguished relation. Lemma~\ref{lem:silting-structure}
gives all structural finiteness conditions and enough injectives.
There are finitely many indecomposable projectives because they
are the indecomposable summands of $P$. The object
$F=P\oplus P[1]$ extension-generates $\mathcal B$ directly from
its definition.

Choose $c\geq1$ such that $\Hom(P,P[-j])=0$ for $j\geq c$.
For $a,d\in\{0,1\}$ and $j\geq c+1$,
\[
 \Hom(P[a],P[d-j])=\Hom(P,P[-(j+a-d)])=0,
\]
since $j+a-d\geq c$. For fixed such $j$ and $a$, applying the
ambient long exact Hom sequence to a defining triangle for $V$
gives $\Hom(P[a],V[-j])=0$. Applying it in the first variable to
a defining triangle for $U$ then gives
$\Hom(U,V[-j])=0$ for all $U,V\in\mathcal B$. Thus ambient
negative Hom gives Definition~\ref{def:negative} with $b=c+1$.

Take the direct sum of the triangles $P\to0\to P[1]$ and
$P[1]\xrightarrow{1}P[1]\to0$. Its first three terms give
\[
 \epsilon:(P\oplus P[1])\longrightarrow P[1]\longrightarrow P[1]
\]
in $\mathcal B$. The middle term is injective by
Lemma~\ref{lem:silting-structure}, and
\[
 r_\epsilon=[P\oplus P[1]]-[P[1]]+[P[1]]=[P]+[P[1]].
\]
Theorem~\ref{thm:general} now shows that (4) implies (1) and that
(1) is equivalent to (2). Condition (2) implies (3) by
Theorem~\ref{thm:ZZ}. Finally (3) implies (4), since the displayed
extriangle shows $[P]+[P[1]]\in\Ex(\mathcal B)$. Hence all four
conditions are equivalent.
\end{proof}

For example, take $\mathcal T=K^b(\operatorname{proj}A)$ and $P=A$.
All the assumptions hold, since the shifted Hom spaces of $A$
vanish in nonzero degrees. In this case $\mathcal B$ consists of
two-term projective complexes up to homotopy and is generally
smaller than $\C_2(A)$.

\section{Examples}\label{sec:examples}
In this section, we give some examples to illustrate the preceding results.

\subsection{A non-hereditary Nakayama example}\label{sec:nakayama}
Let $A=k(1\xrightarrow\alpha2\xrightarrow\beta3)/(\alpha\beta)$,
and let $m=2$. Set $P_i=e_iA$, so $P_3=S_3$,
$\rad P_1=S_2$, and $\rad P_2=S_3$. The algebra is non-hereditary:
$S_1$ has projective resolution
\[
 0\longrightarrow P_3\longrightarrow P_2\longrightarrow P_1
 \longrightarrow S_1\longrightarrow0,
 \qquad \Ext_A^2(S_1,S_3)=k.
\]
Its five indecomposable modules are $S_1,S_2,S_3,P_1,P_2$.
In addition to their degree-zero and degree-one shifts, the category
contains the indecomposable complex
\[
 X=(P_2\longrightarrow P_1),\qquad
 H^{-1}(X)=S_3,\quad H^0(X)=S_1,
\]
where the displayed terms are in degrees $-1,0$ and the differential is
the nonzero radical map.

These are all the indecomposables. Indeed, for $Y\in\C_2(A)$ its
truncation triangle is determined by a connecting morphism
$H^0(Y)\to H^{-1}(Y)[2]$. Among pairs of the five indecomposable
modules, the only nonzero second extension is
$\Ext_A^2(S_1,S_3)=k$. After decomposing the cohomology modules, the
connecting morphism is therefore a matrix over $k$ between copies of
$S_1$ and $S_3[2]$. Changes of basis reduce it to a direct sum of identity
and zero blocks. The identity block gives $X$, and zero blocks give stalk
modules and their shifts. Finally a chain endomorphism of
$P_2\to P_1$ is a pair of equal scalars; there is no nonzero homotopy
$P_1\to P_2$. Thus $\End_{\D}(X)=k$, proving indecomposability.

Order the eleven objects as
\begin{equation}\label{eq:nakorder}
 (Y_1,\ldots,Y_{11})=
 (S_3,P_2,S_2,P_1,S_1,X,S_3[1],P_2[1],S_2[1],P_1[1],S_1[1]).
\end{equation}
Figure~\ref{fig:nakquiver} gives the complete AR quiver; all arrows have multiplicity one.
\begin{figure}[htbp]
\centering
\begin{tikzpicture}[x=1.08cm,y=0.8cm,>=Stealth, every node/.style={font=\small}]
\node (a) at (0,0) {$S_3$};
\node (b) at (1.25,0) {$P_2$};
\node (c) at (2.5,0) {$S_2$};
\node (d) at (3.75,1) {$P_1$};
\node (g) at (3.75,-1) {$S_3[1]$};
\node (f) at (5,0) {$X$};
\node (e) at (6.25,1) {$S_1$};
\node (h) at (6.25,-1) {$P_2[1]$};
\node (i) at (7.5,0) {$S_2[1]$};
\node (j) at (8.85,0) {$P_1[1]$};
\node (k) at (10.2,0) {$S_1[1]$};
\draw[->] (a)--(b); \draw[->] (b)--(c);
\draw[->] (c)--(d); \draw[->] (c)--(g);
\draw[->] (d)--(f); \draw[->] (g)--(f);
\draw[->] (f)--(e); \draw[->] (f)--(h);
\draw[->] (e)--(i); \draw[->] (h)--(i);
\draw[->] (i)--(j); \draw[->] (j)--(k);
\end{tikzpicture}
\caption{The AR quiver of $\C_2(k(1\to2\to3)/\rad^2)$.}\label{fig:nakquiver}
\end{figure}

For reproducibility, we give the morphism calculation establishing this
quiver. Write $H_{ij}=\dimk\Hom(Y_i,Y_j)$ and
$L_{ij}=\dimk\Hom(Y_i,Y_j[-1])$. All their nonzero entries equal one;
Table~\ref{tab:hom} gives their column supports in each row.
\par
\begin{table}[htbp]
\centering
\caption{Nonzero entries of the Hom and negative Hom matrices.}\label{tab:hom}
\begin{tabular}{ccl}
\toprule
$i$ & $\{j:H_{ij}=1\}$ & $\{j:L_{ij}=1\}$\\
\midrule
1 & $1,2$ & $6,7,8$\\
2 & $2,3,4$ & $8,9,10$\\
3 & $3,4,6,7$ & $9,10$\\
4 & $4,5,6$ & $10,11$\\
5 & $5,9$ & $11$\\
6 & $5,6,8,9$ & $11$\\
7 & $6,7,8$ & $\varnothing$\\
8 & $8,9,10$ & $\varnothing$\\
9 & $9,10$ & $\varnothing$\\
10 & $10,11$ & $\varnothing$\\
11 & $11$ & $\varnothing$\\
\bottomrule
\end{tabular}
\end{table}
\par
To obtain the table over any field, use the projective complexes
$(P_3\to P_2\to P_1)$ for $S_1$, $(P_3\to P_2)$ for $S_2$,
and stalk projectives for $S_3,P_1,P_2$, together with their shifts and
the displayed complex $X$. The rightmost term is in degree zero.
The spaces $\Hom_A(P_i,P_j)$ are one-dimensional for $i=j$ or $i=j+1$
and zero otherwise. The composite of two radical maps is zero.
For complexes $B,C$, form
\[
 \mathcal H^q(B,C)=\bigoplus_t\Hom_A(B^t,C^{t+q}),\qquad
 d(f)=d_Cf-(-1)^q f d_B.
\]
Then $H^0\mathcal H$ and $H^{-1}\mathcal H$ give the two tables.
Each nonzero class is represented by a single scalar or by equal scalars
on adjacent terms. Modulo homotopies, the only nonzero spaces
$\rad^2(Y_i,Y_j)$ occur at
\[
 (i,j)=(2,4),(3,6),(4,5),(6,9),(7,8),(8,10).
\]
For example, the two paths $S_2\to P_1\to X$ and
$S_2\to S_3[1]\to X$ give the same one-dimensional space up to sign;
the two paths from $X$ to $S_2[1]$ behave likewise. The other four pairs
have a unique nonzero two-arrow composite. All remaining products of
two radical classes vanish, and $\rad^3=0$.
Thus $\rad/\rad^2$ gives exactly the arrows drawn above.

The eight AR extriangles and their relations are recorded in Table~\ref{tab:meshes}.
The projective objects are $Y_1,Y_2,Y_4$.
\par
\begin{table}[htbp]
\centering
\caption{Auslander--Reiten extriangles for the eleven-object example.}\label{tab:meshes}
\begin{tabular}{cccc}
\toprule
Relation & First term & Middle term & Endpoint\\
\midrule
$r_1$ & $S_3$ & $P_2$ & $S_2$\\
$r_2$ & $S_3[1]$ & $X$ & $S_1$\\
$r_3$ & $S_2$ & $P_1\oplus S_3[1]$ & $X$\\
$r_4$ & $P_2$ & $S_2$ & $S_3[1]$\\
$r_5$ & $P_1$ & $X$ & $P_2[1]$\\
$r_6$ & $X$ & $S_1\oplus P_2[1]$ & $S_2[1]$\\
$r_7$ & $S_1$ & $S_2[1]$ & $P_1[1]$\\
$r_8$ & $S_2[1]$ & $P_1[1]$ & $S_1[1]$\\
\bottomrule
\end{tabular}
\end{table}
\par
Here each relation means first minus middle plus endpoint. Direct
integer elimination gives $\det(H-L)=1$. By
Lemma~\ref{lem:ARsupport}, the relation of an AR extriangle ending
at $Y_j$ is the unique vector $r$ satisfying $(H-L)r=e_j$.
Each displayed relation satisfies this equation, so the listed
relations are exactly the AR relations. The irreducible maps in
the quiver determine the middle terms of the minimal right almost
split maps; the relation then determines each first term.

In simple coordinates $(S_1,S_2,S_3)$ the Euler matrix is
\[
 E=\begin{pmatrix}
 0&0&0&1&1&1&0&0&0&-1&-1\\
 0&1&1&1&0&0&0&-1&-1&-1&0\\
 1&1&0&0&0&-1&-1&-1&0&0&0
 \end{pmatrix}.
\]
Let $R$ be the $11\times8$ matrix with columns $r_1,\ldots,r_8$ in
the order \eqref{eq:nakorder}. Then $ER=0$. Appending the columns
$[S_1],[S_2],[S_3]$ gives an integer matrix of determinant one.
Since $E$ sends those last three columns to the identity matrix,
\[
 0\longrightarrow\mathbb Z^8\xrightarrow{R}\mathbb Z^{11}
 \xrightarrow{E}\mathbb Z^3\longrightarrow0
\]
is exact over $\mathbb Z$. In particular the eight AR relations form an
integral basis of $\Ex(\C_2(A))$.

Finally choose $J=S_1\oplus P_1\oplus P_2$ as the injective envelope
of $S=S_1\oplus S_2\oplus S_3$; its quotient is $Q=S_1\oplus S_2$.
The relation $r_\epsilon$ becomes
\begin{align*}
 r_\epsilon={}&[S_3]+[S_2]+[S_1]+2[S_3[1]]-[P_2[1]]\\
             &\quad+3[S_2[1]]-[P_1[1]]+2[S_1[1]]\\
 ={}&r_1+2r_2+r_3+r_4+r_5+2r_6+r_7+2r_8.
\end{align*}
All residual division rings equal $k$. Evaluating with $H-L$ at the
eight endpoints gives $(1,2,1,1,1,2,1,2)$, exactly the defect coefficients
in \eqref{eq:ARcoefficients}.

\subsection{Counts for radical-square-zero linear Nakayama algebras}\label{sec:counts}
Let $A_n=k(1\to\cdots\to n)/\rad^2$, $n\geq2$.
These are piecewise hereditary of Dynkin type, within the classification
of Happel--Seidel \cite{HS} recalled in \cite[Section 4.1]{Rundsveen}.
Section~\ref{sec:nakayama} is the
case $n=3$, and the formula below gives $N(3,2)=11$.

\begin{proposition}\label{prop:nakcount}
For every $m\geq1$,
\begin{equation}\label{eq:nakcount}
 N(n,m):=\#\ind\C_m(A_n)
 =m(2n-1)+\sum_{d=1}^{\min(m-1,n-2)}(n-1-d)(m-d).
\end{equation}
Consequently $\operatorname{rank}\Ex(\C_m(A_n))=N(n,m)-n$.
For $m\geq n-1$ this simplifies to
$N(n,m)=m\binom{n+1}{2}-\binom n3$.
\end{proposition}
\begin{proof}
The simple modules have finite projective resolutions of length at most
$n-1$, so every bounded derived object has a bounded projective
representative. Remove contractible summands to make it minimal.
Since the only radical maps between projectives are
$P_i\to P_{i-1}$, the coefficient vector spaces in a minimal complex
separate into diagonals indexed by the sum of the degree and the vertex
index. On each diagonal they form a representation of a linearly oriented
quiver of type $A_n$. The differential-square equation imposes no extra
condition on these coefficient matrices, because the composite of two
radical maps in $A_n$ is already zero. The interval decomposition of
representations of a line, obtained by bases adapted to successive images
and kernels, therefore decomposes the complex into shifts of
\[
 P_b\longrightarrow P_{b-1}\longrightarrow\cdots\longrightarrow P_a,
 \qquad 1\leq a\leq b\leq n.
\]
Every such interval has scalar endomorphisms and no degree-minus-one
endomorphism giving a nonzero homotopy, hence is indecomposable. Minimal
bounded complexes are homotopy equivalent only if they are isomorphic:
a homotopy equivalence is invertible modulo the radical in each degree,
so is degreewise invertible. Thus the intervals and their shifts give
distinct indecomposables.

The $n$ intervals of length one are projective stalks and contribute
$mn$ objects. For $a<b=n$, an interval has only the top cohomology
$S_a$ and contributes $m$ shifts, giving $m(n-1)$ more objects.
If $a<b<n$ and $d=b-a$, its cohomology consists of $S_a$ in the top
degree and $S_{b+1}$ exactly $d$ degrees lower. Both lie in the cohomological interval
for precisely $m-d$ shifts when $d<m$, and for none otherwise.
There are $n-1-d$ such intervals. Summing gives \eqref{eq:nakcount}.
The rank statement follows from Proposition~\ref{prop:kernel}; evaluating
the finite sums gives the last formula. Its leading coefficient
$\binom{n+1}{2}$ is the number of indecomposable modules of Dynkin
type $A_n$, consistent with the interval description above.
\end{proof}

Table~\ref{tab:counts} lists $N(n,m)$ with the relation rank in parentheses.
\par
\begin{table}[!htbp]
\centering
\caption{Indecomposable counts and relation ranks for the radical-square-zero family.}\label{tab:counts}
\begin{tabular}{ccccc}
\toprule
$n\backslash m$ & 1 & 2 & 3 & 4\\
\midrule
3 & $5\ (2)$ & $11\ (8)$ & $17\ (14)$ & $23\ (20)$\\
4 & $7\ (3)$ & $16\ (12)$ & $26\ (22)$ & $36\ (32)$\\
5 & $9\ (4)$ & $21\ (16)$ & $35\ (30)$ & $50\ (45)$\\
\bottomrule
\end{tabular}
\end{table}
\par

\FloatBarrier

\subsection{The Kronecker algebra and a relation outside \texorpdfstring{$\AR$}{AR}}
Let $A=kQ$, where $Q$ has two arrows $\alpha,\beta:1\to2$.
For $t\geq1$, consider the representation
\[
 M_t=(k^t,k^{t+1};\alpha_t,\beta_t),\qquad
 \alpha_t(x)=(x,0),\quad\beta_t(x)=(0,x).
\]
These representations are pairwise nonisomorphic indecomposables over
every field. Indeed, an endomorphism $(a,b)$ satisfies
$b\alpha_t=\alpha_t a$ and $b\beta_t=\beta_t a$. Comparing the
overlapping columns gives
\[
 a_{1j}=0\ (j>1),\quad a_{tj}=0\ (j<t),\quad
 a_{rj}=a_{r-1,j-1}\ (r,j\geq2).
\]
These equations force $a$ to be a scalar matrix, and the two arrow
images span $k^{t+1}$, forcing the same scalar for $b$. Hence
$\End(M_t)=k$. Their different dimension vectors distinguish them.
Since $\Hom(A,M_t)\ne0$ for all $t$, $\C_m(A)$ is not locally finite
for any $m$. Nevertheless, Proposition~\ref{thm:K0} gives
$K_0(\C_m(A))\simeq\mathbb Z^2$ for every $m$.

One can exhibit an explicit relation which is not generated by AR
extriangles. Let
$I_2=(k^2,k;\pi_1,\pi_2)$, with the two coordinate projections.
It is the indecomposable injective at vertex 2, and
\[
 0\longrightarrow S_2\longrightarrow I_2
   \longrightarrow S_1^{\oplus2}\longrightarrow0
\]
is exact. After shifting by $m-1$, its relation is
\[
 r_m=[S_2[m-1]]-[I_2[m-1]]+2[S_1[m-1]]\in\Ex(\C_m(A)).
\]
Put $U_t=M_t[m-1]$. All negative Hom terms with $j>0$ in
$\chi^-_{U_t}$ vanish on these three shifted modules. Also,
\[
 \dimk\Hom(M_t,S_2)=0,\qquad
 \dimk\Hom(M_t,I_2)=t+1,\qquad
 \dimk\Hom(M_t,S_1)=t.
\]
For the first equality, a map at vertex 2 must annihilate both arrow
images, whose sum is the whole space. The last equality comes from an
arbitrary functional on $k^t$. The middle equality follows either by
the same matrix calculation or from
$\Hom(M_t,D(Ae_2))\simeq D(M_te_2)$.
Consequently,
\[
 \chi^-_{U_t}(r_m)=0-(t+1)+2t=t-1>0\qquad(t\geq2).
\]
Every finite integral combination of AR relations has finite support
under these evaluations by Lemma~\ref{lem:ARsupport}. Since the $U_t$ are
pairwise nonisomorphic, it follows that $r_m\notin\AR(\C_m(A))$.
More strongly, if $q\ne0$ is an integer, then
$\chi^-_{U_t}(q r_m)=q(t-1)$ still has infinite support, so
$q r_m\notin\AR(\C_m(A))$. Hence
\[
 \mathbb Z\hookrightarrow\Ex(\C_m(A))/\AR(\C_m(A)),\qquad
 q\longmapsto q\,\overline{r_m}
\]
is injective.

\section*{Acknowledgements}

The author is grateful to Professor Bin Zhu for suggesting the
question of whether an extended module category is finit if and only if its Grothendieck relation
subgroup is generated by Auslander--Reiten extriangles.
The author also sincerely thanks him for his helpful discussions,
valuable suggestions, and generous support throughout the
research of this question.


\PrintDeclaration{AI use statement}{\AIUseStatement}


\begin{thebibliography}{99}
\raggedright
\interlinepenalty=10000
\bibitem{AET}
T. Adachi, H. Enomoto and M. Tsukamoto, \emph{Intervals of $s$-torsion
pairs in extriangulated categories with negative first extensions},
Mathematical Proceedings of the Cambridge Philosophical Society
\textbf{174} (2023), 451--469.

\bibitem{AI}
T. Aihara and O. Iyama, \emph{Silting mutation in triangulated categories},
Journal of the London Mathematical Society \textbf{85} (2012), 633--668.

\bibitem{Auslander}
M. Auslander, \emph{Relations for Grothendieck groups of Artin algebras},
Proceedings of the American Mathematical Society \textbf{91} (1984),
336--340.
\bibitem{Butler}
M. C. R. Butler, \emph{Grothendieck groups and almost split sequences},
in \emph{Integral Representations and Applications}, Lecture Notes in
Mathematics \textbf{882}, Springer, 1981, 357--368.
\bibitem{Enomoto}
H. Enomoto, \emph{Relations for Grothendieck groups and
representation-finiteness}, Journal of Algebra \textbf{539} (2019),
152--176.

\bibitem{FedeleAR}
F. Fedele, \emph{Auslander--Reiten $(d+2)$-angles in subcategories and a
$(d+2)$-angulated generalisation of a theorem by Br\"uning}, Journal
of Pure and Applied Algebra \textbf{223} (2019), 3554--3580.

\bibitem{FedeleK}
F. Fedele, \emph{Grothendieck groups of triangulated categories via
cluster tilting subcategories}, Nagoya Mathematical Journal
\textbf{244} (2021), 204--231.

\bibitem{GNP}
M. Gorsky, H. Nakaoka and Y. Palu, \emph{Positive and negative extensions
in extriangulated categories}, Mathematische Zeitschrift \textbf{313}
(2026), Article 89.

\bibitem{Gupta}
E. Gupta, \emph{On $d$-term silting objects, torsion classes, and cotorsion
classes}, \href{https://arxiv.org/abs/2407.10562v3}{arXiv:2407.10562v3}, 2026.
\bibitem{GZ}
E. Gupta and Y. Zhou, \emph{Semibricks and wide subcategories in extended
module categories}, \href{https://arxiv.org/abs/2511.08157v1}{arXiv:2511.08157v1}, 2025.
\bibitem{HS}
D. Happel and U. Seidel, \emph{Piecewise hereditary Nakayama algebras},
Algebras and Representation Theory \textbf{13} (2010), 693--704.

\bibitem{Haugland}
J. Haugland, \emph{Auslander--Reiten triangles and Grothendieck groups of
triangulated categories}, Algebras and Representation Theory
\textbf{25} (2022), 1379--1387.

\bibitem{HLZ}
J. He, Y.-Z. Liu and P. Zhou, \emph{Auslander--Reiten $(n+2)$-angles
and local finiteness},
\href{https://arxiv.org/abs/2608.26653v1}{arXiv:2608.26653v1}, 2026.
\bibitem{INP}
O. Iyama, H. Nakaoka and Y. Palu, \emph{Auslander--Reiten theory in
extriangulated categories}, Transactions of the American Mathematical
Society, Series B \textbf{11} (2024), 248--305.
\bibitem{MNO}
N. Mochizuki, H. Nakaoka and Y. Ogawa, \emph{Extended heart construction
(I): The heart of $n$-cotorsion pairs on triangulated categories},
\href{https://arxiv.org/abs/2608.03212v1}{arXiv:2608.03212v1}, 2026.
\bibitem{NP}
H. Nakaoka and Y. Palu, \emph{Extriangulated categories, Hovey twin
cotorsion pairs and model structures}, Cahiers de Topologie et G\'eom\'etrie
Diff\'erentielle Cat\'egoriques \textbf{60} (2019), 117--193.
\bibitem{PPPP}
A. Padrol, Y. Palu, V. Pilaud and P.-G. Plamondon,
\emph{Associahedra for finite type cluster algebras and minimal relations
between $g$-vectors}, Proceedings of the London Mathematical Society
\textbf{127} (2023), 513--588.

\bibitem{Rundsveen}
E. S. Rundsveen, \emph{A note on the $m$-extended module categories of
Nakayama algebras},
\href{https://arxiv.org/abs/2601.14843v1}{arXiv:2601.14843v1}, 2026.
\bibitem{Stacks}
The Stacks Project Authors, \emph{The Stacks Project}, Derived Categories,
Section 13.28, Lemma 13.28.2 (Tag 0FCP), 2026.
\url{https://stacks.math.columbia.edu/tag/0FCP}.
\bibitem{WZ}
J. Wei and Y. Zhou, \emph{AIR tilting subcategories of extended hearts},
\href{https://arxiv.org/abs/2512.13218v3}{arXiv:2512.13218v3}, 2026.
\bibitem{XZrelations}
J. Xiao and B. Zhu, \emph{Relations for the Grothendieck groups of
triangulated categories}, Journal of Algebra \textbf{257} (2002), 37--50.

\bibitem{XZlocal}
J. Xiao and B. Zhu, \emph{Locally finite triangulated categories},
Journal of Algebra \textbf{290} (2005), 473--490.
\bibitem{PZhou}
P. Zhou, \emph{Grothendieck groups and Auslander--Reiten $(d+2)$-angles},
Journal of Algebra \textbf{568} (2021), 512--528;
\href{https://arxiv.org/abs/1912.11397v2}{arXiv:1912.11397v2}.

\bibitem{ZhouZhu}
P. Zhou and B. Zhu, \emph{Auslander–Reiten relations and local finiteness in extriangulated categories} 2026

\bibitem{Zhou}
Y. Zhou, \emph{Tilting theory for extended module categories},
\href{https://arxiv.org/abs/2411.15473v2}{arXiv:2411.15473v2}, 2025.

\bibitem{ZZ}
B. Zhu and X. Zhuang, \emph{Grothendieck groups in extriangulated
categories}, Journal of Algebra \textbf{574} (2021), 206--232.


\end{thebibliography}
\end{document}